\documentclass{amsart}

\usepackage{amssymb}
\usepackage{graphicx}

\usepackage[margin=1in]{geometry}  % set the margins to 1in on all sides
\usepackage{graphicx}              % to include figures
\usepackage{amsmath}               % great math stuff
\usepackage{amsfonts}              % for blackboard bold, etc
\usepackage{amsthm}                % better theorem environments

\begin{document}

\nocite{*}
\title{Spot it!$^{\textregistered}$  Solitaire}

\author{Donna A. Dietz\\
Department of Mathematics and Statistics\\
 American University\\ Washington, DC, USA}
\email{donna.dietz@american.edu}

\maketitle

\begin{abstract}The game of Spot it$^{\textregistered}$ is based on an order 7 finite projective plane. This article presents a solitaire challenge: extract an order 7 affine plane and arrange those 49 cards into a square such that the symmetries of the affine and projective planes are obvious. The objective is not to simply create such a deck already in this solved position. Rather, it is to solve the inverse problem of arranging the cards of such a deck which has already been created  shuffled.
\end{abstract}

\section {Introduction}

As I was shopping online one day, an advertisement for Spot
it!$^{\textregistered}$ 
caught my eye. This game is played with 55
circular cards, each card having several images, and each pair of
cards having exactly one common image. Several games can be played
with the deck, all involving multiple players trying to be the first
to spot the matching image between two cards.\\

\begin{figure}
\begin{center}
\includegraphics[width=3in]{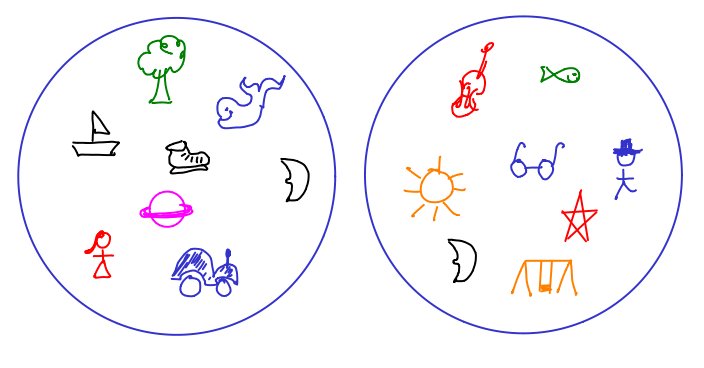}\\
\end{center}
\caption{These two cards have an image in common.}
\label{cardDemo}
\end{figure}

Naturally, I purchased the deck, which is made by {\it Blue Orange Games},  and quickly discovered that the deck is two cards shy of fully representing an order 7 finite projective plane. An order 7 finite projective plane has 57 points, 57 lines, 8 lines per point and 8 points per line. Each pair of lines is incident with a unique point, and each pair of points is incident with a unique line. The cards represent the points and the images represent the lines, or vice versa due to duality. Those readers who are unfamiliar with finite geometry models may wish to read Maxime Bourrigan's online article about the finite geometry behind Spot it$^{\textregistered}$ \cite{bourrigan}. (Due to the nature of the article, a Google translation suffices.)
For those whose interest in this topic is beyond the content of this discussion, I recommend \cite{dorwart} for students who have not yet mastered abstract algebra, and \cite{hirschfeld} and \cite{hughes} for those who have.
\\ 

It seemed a natural course of action to create the two missing cards and then proceed to arrange the cards into a configuration which would make it easy to demonstrate the order 7 finite projective plane. I didn't realize how fun and challenging this would be. I'm hoping the ``rules'' (and solution) of this single-player challenge will be entertaining to mathematicians and game-lovers alike. \\

This discussion is arranged so that those wishing fewer clues can read fewer sections, thus leaving more of the fun for themselves.\\

\section{The rules of the game}

To play, you need a deck of Spot it$^{\textregistered}$ cards or a similar set of cards which represent a finite projective plane of order 7 (or some other prime order, $n$). Begin by removing all cards having a specific image. The remaining cards form an affine plane of the order 7 (or $n$).  (In my Spot it$^{\textregistered}$ deck, the two missing cards both contain a snowman.  So, if I simply remove all the snowman cards at this step, I do not need to actually find the two missing cards in order to continue.)  If the two missing cards are replaced, cards of any image may be removed at this step, thanks to the self-symmetry of the finite projective plane.\\

The ultimate goal is to lay out the  49 ($n^2$) remaining cards in a 7 by 7 ($n$ by $n$) grid so that this one rule is satisfied: Pick any two cards in the grid. Let their positions in the grid be $(x,y)$ and $(x+h,y+k)$ with $x$ and $y$ numbered between 0 and 6 inclusive. The common image found on these two cards must also be at position $(x+2h \mod{7}, y+2k \mod{7})$  (or $(x+2h \mod{n}, y+2k \mod{n})$).\\

\begin{figure}
\begin{center}
\includegraphics[width=2in]{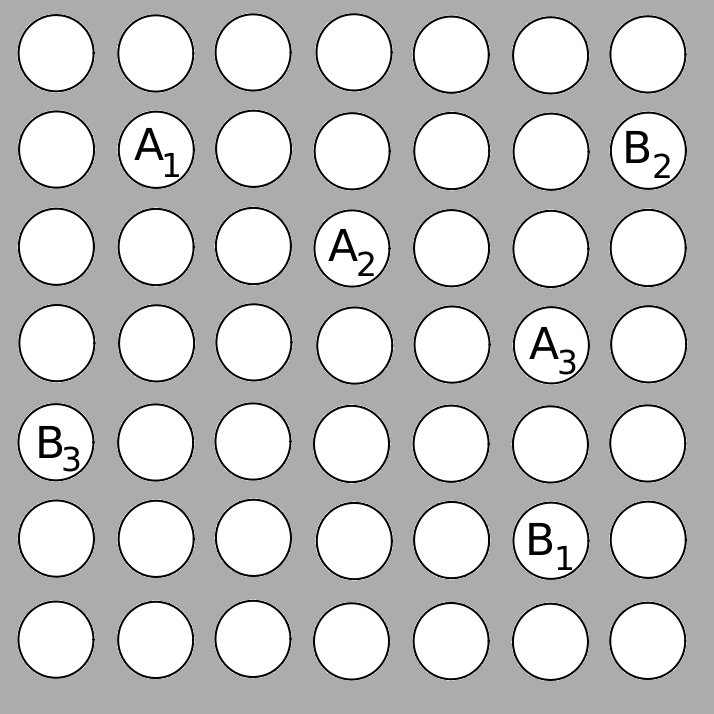}\\
\end{center}
\caption{Demonstration of two collections of three images.}
\label{gridDemo}
\end{figure}

For example, in Figure \ref{gridDemo}, the circles represent cards. The image which is common between the cards in positions $A_1$ and $A_2$ must also be present on the card in position $A_3$. (Likewise for the cards in positions $B_1$, $B_2$, and $B_3$.) Clearly, since 7 ($n$) is prime, and therefore all elements are generators in $\mathbb{Z}_n$ (or $\mathbb{Z}_n$), 
there will be 7 ($n$) such images in a set. This also implies that each row and each column will have a common image.\\

\begin{figure}
\begin{center}
\includegraphics[width=4in]{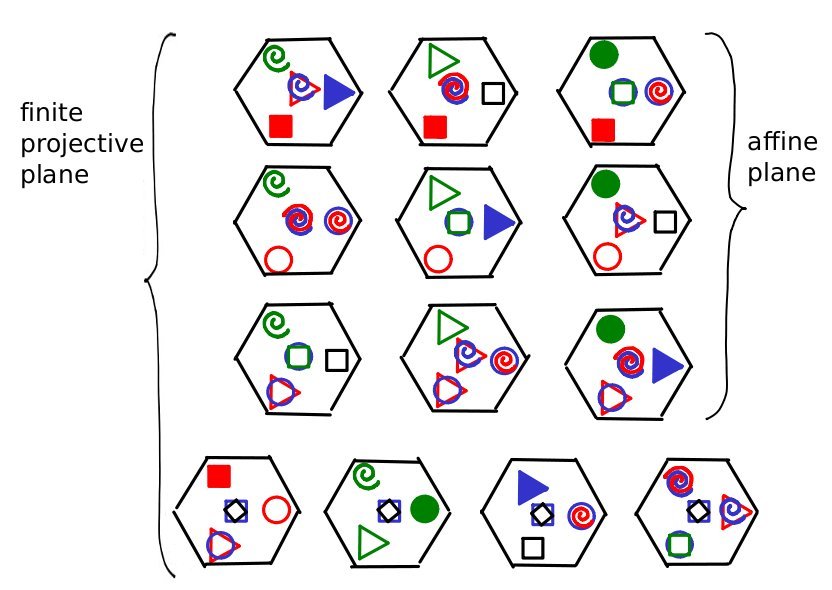}\\
\end{center}
\caption{The solved challenge for the $n=3$ case.}
\label{n3case}
\end{figure}

In Figure \ref{n3case}, the solved $n=3$ case is shown. Look at each of the 13 symbols, one by one, to see how the lines fit together.\\

The technique for creating such a deck of cards should be obvious
now. (Simply use the one rule to generate common symbols on the grid
as specified. Then, symbols sharing a slope cannot coexist on a card
in the grid, but they will be found together on the cards outside
the grid. All the cards outside the grid will also contain one additional
common symbol.)  However, the goal here is not to create such a deck
but rather to properly arrange an already existing deck!\\

Those readers desiring the maximal fun should now attempt to solve this inverse problem without reading further!  One more warning will be given after the easy clues are presented.\\

\section{Finding the missing cards}
The techniques presented here may be generalized for more than two missing cards, but I will presume the reader has a Spot it$^{\textregistered}$ deck which is missing two cards. First, you will need a list of all the images in your deck. The easiest way to do this is to locate an image which occurs 8 times (that is, $n+1$), and pull out all those cards.  For example, if you have 8 spiders, you can pull out the spider ``set''. Since each of those cards has a spider in common, there can be no other commonalities in any of those 8 cards. Thus, there must be $8 \times 7 + 1$ or all 57 images on those 8 cards. \\

Next, make an accounting of the frequencies for each image. There should be one image which is present only 6 times, while 14 images are present 7 times. All other images should be present 8 times, making complete sets. The image which is missing twice must occur on both missing cards. 
Without loss of generality, place one additional missing image on one of the cards you are creating. For example, a maple leaf. Then, for the remaining 13 images, search the entire deck to see if that image occurs with the maple leaf or not. If it does, it cannot do so again, so it must go on the other new card. If, however, it does not occur with the maple leaf, it should!  So, it should go on the new card containing the maple leaf!

\section{The initial setup of the grid}
As was mentioned in the rules, you need to first pull out a set of 8 cards which share an image. (Or, if you have not created the missing cards, you should pull out 6 cards sharing the twice-missing image.) We will refer to this image as ``infinity'', because its set will be the set which differentiates between our affine plane and our finite projective plane! The grid we are building is the respresentation of the affine plane associated with this projective plane.\\

The cards in the infinity set each contain two things. They contain the infinity image, and they contain a collection of parallel images of the resulting affine plane. In other words, any images from the affine plane (grid) which never occur together will definitely occur together on an infinity card. \\

Select two such infinity cards, one to keep track of the rows of your grid, and another to keep track of the columns of your grid. Arrange your grid now so that each column contains a common image and each row contains a common image.
Note that there are $57 \times 8 \times 7 \times 7! \times 7!$ ways to do this. There are 57 images to pick as the infinity image, then 8 cards from that infinity set which can be used to define the rows. Once the rows are chosen, 7 cards remain to define the columns. Once the rows and columns are determined, there are 7! ways to order  rows and 7! ways to order columns. Now, by switching rows and/or columns, place a common symbol on just the main diagonal of the grid. All ``moves'' henceforth will consist of swapping two rows or swapping two columns (or any equivalent). We know from abstract algebra that all permutations can be formed by swaps, so swapping rows/columns is sufficient, but there are often faster techniques in practice. This will be obvious to the reader who attempts this.\\

The next objective is to get a common symbol on just the main counter-diagonal. As the final easy clue, you now have available to you the option to do as you wish to the rows so long as you do exactly the same things to the columns of the same indices. This will allow you to maintain your main diagonal while still giving you enough freedom to finish the puzzle. We will, without loss of generality, freeze the middle card and do not move the middle row or middle column. (We are on a torus anyhow, and freezing the middle card gives us the advantage of symmetry which will become apparent as the solutions progresses.)  This gives $6!=720$ remaining grid arrangements, 6 of which are valid solutions. At this point, if you wish to have any fun with the puzzle, you should stop reading. This is your last warning.\\

\section{Finding the counterdiagonal and finishing the puzzle}

The somewhat surprising result is that once you have fixed the image on your main diagonal and have frozen the middle card, your counterdiagonal image is already determined! It has to be one of the 8 images on your middle card, but it cannot be the image of that row or of that column or of the diagonal. So, there appear to be 5 options remaining. But that is not so! Only one will work, and any attempt to set one of the other 4 images as the counterdiagonal will end in frustration, whereas an attempt to set the correct image as the counterdiagonal will end smoothly and quickly with success.  So, how do we figure out which image will work?\\

\begin{figure}
\begin{center}
\includegraphics[width=3in]{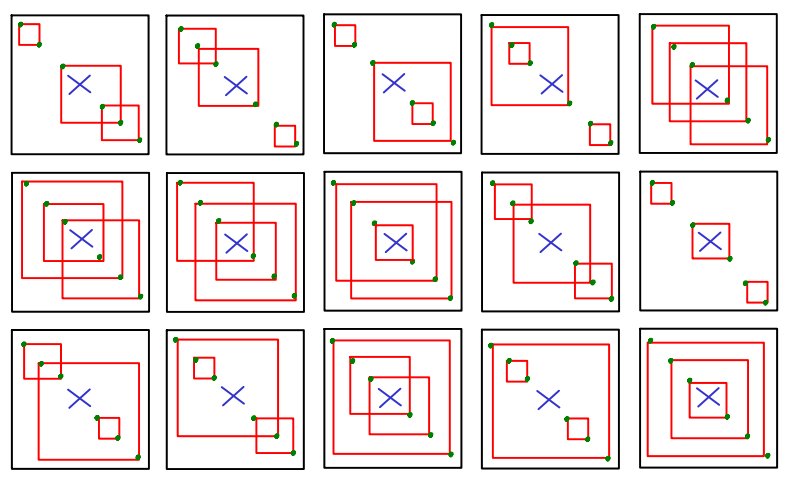}\\
\end{center}
\caption{The 15 ways the three squares may interlace.}
\label{fifteenSq}
\end{figure}

Let us imagine the reverse. What happens when we start with the
correct arrangement of cards and then scramble them so that
corresponding rows/columns are always swapped at the same time (one
after the other, in either order)?  You can define concentric squares
around the middle element which live on the main diagonal of the
grid. The corners of the squares are the elements of the
counterdiagonal and the diagonal. Strangely enough, the (group)
actions of swapping the same rows as columns actually maintains these
three squares along the main diagonal, although they are not as
concentric squares except in the case of the solved counterdiagonal!
So, in order to determine which of the 5 candidate images should be
used as the counterdiagonal image, simply start at the first element
and count down (or across) until the candidate symbol is located. The
two counts should be the same when you have found the correct
image. This count must always be equal, starting from any point on the
diagonal, counting over and down (or back and up). What this tells us
is that there is a natural pairing of our 6 counter-diagonal elements
which can be detected when the correct candidate image is
picked. Combinatorially, 6 elements in 3 pairs can occur 15 ways. This
is demonstrated in Figure \ref{fifteenSq} which shows how these
pairings can be manifest in the grid.\\

When you get to this point in the puzzle, it is a beautiful thing to find these squares, reassuring you that the solution is near! Simply swap rows and equivalent columns until the counterdiagonal is set!\\

\section{Final moves}

Once the counterdiagonal is set, remaining  moves must not only have paired row and column moves, but must also maintain symmetry between right and left (as well as up and down). For example, if you swap rows/columns 3 and 5, this is already balanced and will not disrupt the counterdiagonal. However, if you wish to swap columns 1 and 2, you must also swap columns 6 and 7 (as well as rows 1 and 2 along with rows 6 and 7). By tightening these orbits, you are closing in on one of the 6 final solutions. \\

Recall that there were $6!=720$ ways to arrange the grid of cards while maintaining the diagonal once the middle card was fixed.  This count agrees with another counting scheme. There are 15 ways for the counterdiagonal squares to be arranged, only one of which is useful for a solution, of course. But once you restrict yourself to that pattern of squares, you expect to find 720/15 or 48 arrangements for the cards. This is in fact the case. Even though you have set down the pattern for the squares, you have not defined which square ought to go where, and there are 6 ways to do that. Each square may additionally be in one of two orientations, as it is legal to rotate it 180 degrees. $2^3=8$ and $6\times 8=48$. Of those 48 possible arrangements of the cards, 6 are solutions. For example, you have the freedom to choose any one square's location and orientation, but the rest is then predetermined. For simplicity, we presume the innermost square is set properly, that is, the 9 cards in the middle of the grid are now fixed.\\

\begin{figure}
\begin{center}
\includegraphics[width=1in]{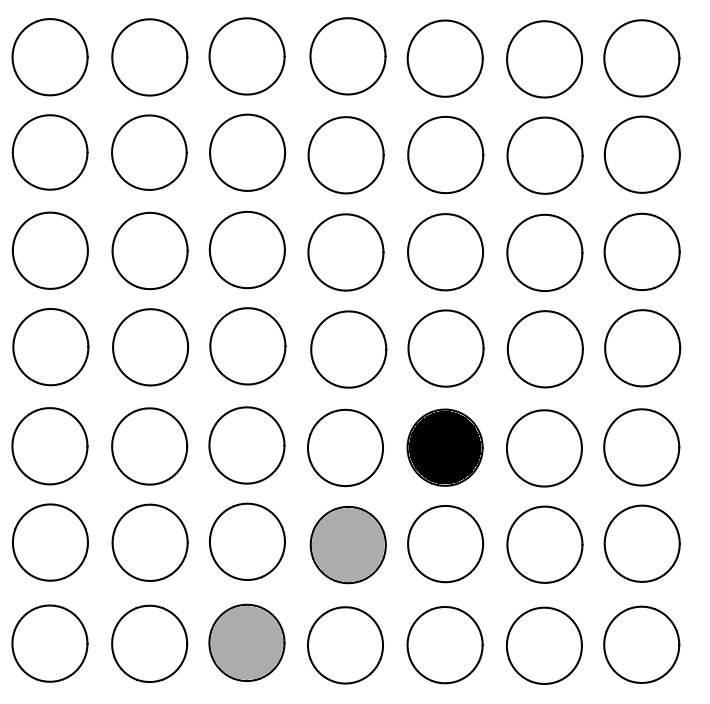}\\
\end{center}
\caption{Propogation of a known image on a fixed square.}
\label{prop}
\end{figure}

Now, it should be easy to slide the cards into their final positions. 
We  use row/column/right/left symmetric moves to assure that some attribute is maintained. For example, we now know which family of images defines the counterdiagonals. As suggested in Figure \ref{prop}, the counterdiagonal image can be propogated from the card indicated by the black circle to nearby cards indicated by the gray circles. In a few moves, you will see before you, a perfectly arranged affine plane.\\


\begin{thebibliography}{9}

\bibitem{bourrigan}
  Maxime Bourrigan,
  \emph{Dobble et la g\'{e}om\'{e}trie finie}.
  Images des Math\'{e}matiques, CNRS, 2011.


\bibitem{dorwart}
  Harold L. Dorwart,
  \emph{The Geometry of Incidence}.
  Prentice-Hall, Inc,
  Englewood Cliffs, N.J, 1966.

\bibitem{hirschfeld}
  J. W. P. Hirschfeld,
  \emph{Projective Geometries over Finite Fields}.
  Second Edition.
  Oxford Science Publications,
  Clarendon Press, Oxford, 1998.

\bibitem{hughes}
  Daniel R. Hughes and Fred C. Piper,
  \emph{Projective Planes}.
  Springer-Verlag, New York, 1973.


\end{thebibliography}
\end{document}